\documentclass [a4,11pt]{amsart}

\usepackage{amsmath,amssymb,amsfonts,enumerate,amsthm, amscd,}

\newtheorem{thm}{Theorem}[section]
\newtheorem{cor}[thm]{Corollary}
\newtheorem{lem}[thm]{Lemma}
\newtheorem{prop}[thm]{Proposition}

\newtheorem{exam}[thm]{Example}

\def\proof{{\parindent0pt {\bf Proof.\ }}}

\theoremstyle{definition}

\theoremstyle{remark}

\theoremstyle{Definition and Notation}

\begin{document}

\bibliographystyle{amsplain}


\title[When every principal ideal is flat]{When every principal ideal is flat}

\author{Fatima Cheniour}
\address{Fatima Cheniour\\Department of Mathematics, Faculty of Science and Technology of Fez, Box 2202, University S.M. Ben Abdellah Fez, Morocco.
$$ E-mail\ address:\ cheniourfatima@yahoo.fr$$}

\author{Najib Mahdou}
\address{Najib Mahdou\\Department of Mathematics, Faculty of Science and Technology of Fez, Box 2202, University S.M. Ben Abdellah Fez, Morocco.
 $$E-mail\  address:\ mahdou@hotmail.com$$}


\begin{abstract}
 This paper deals with well-known notion of $PF$-rings, that is, rings in which principal ideals are flat.
 We give a new characterization of $PF$-rings. Also, we provide a necessary and sufficient condition for $R\bowtie I$
(resp., $R/I$ when $R$ is a Dedekind domain or $I$ is a primary
ideal) to be $PF$-ring. The article includes a brief discussion of
the scope and precision of our results.

\end{abstract}

\maketitle

 \begin{section} {Introduction}
 All rings considered in this paper are assumed to be commutative with identity elements and all modules are
 unitary. \\
 We start by recalling some definitions. \\

 A ring $R$ is called a $PF$-ring if principal ideals of $R$ are
flat. Recall that $R$ is a $PF$-ring if and only if $R_Q$ is a
domain for every prime (resp., maximal) ideal $Q$ of $R$. For
example, any domain and any semihereditary ring is a $PF$-ring
(since a localization of a semihereditary ring by a prime (resp.,
maximal) ideal is a valuation domain). Note that a $PF$-ring is
reduced by \cite[Theorem 4.2.2 , p. 114]{G}. See for
instance \cite{G, H}. \\

 An R-module $M$ is called $P$-flat if, for any $(s,x) \in
R\times M$ such that $ sx=0$, then $ x\in$ $(0:s)M$. If M is flat,
then M is naturally $P$-flat. When $R$ is a  domain, $M$ is
$P$-plat if and only if it is torsion-free. When $R$ is an
arithmetical ring, then any P-flat module is flat (by \cite[p.
236]{C}). Also, every $P$-flat cyclic module is flat (by
\cite[Proposition 1(2)]{C}). See
for instance \cite{C, G}. \\

 The amalgamated duplication of a ring $R$ along an ideal $I$ is a
ring that is defined as the following subring with unit element
$(1,1)$ of $R\times R$:
\begin{eqnarray*}
   &R\bowtie I =\{(r,r+i)/r\in R,i\in I\}.&
\end{eqnarray*}
This construction has been studied, in the general case, and from
the different point of view of pullbacks, by  D'Anna and Fontana
\cite{DF2}. Also, in \cite{DF1}, they  have considered the case of
the amalgamated duplication of a ring, in not necessarily
Noetherian setting, along a multiplicative canonical ideal in the
sense of \cite{HHP}. In \cite{D} D'Anna has studied some
properties of $R\bowtie I $, in order to construct reduced
Gorenstein rings associated to Cohen-Macaulay rings and has
applied this construction to curve singularities. On the other
hand,  Maimani and  Yassemi, in \cite{MY}, have studied the
diameter and girth of the zero- divisor of the ring $R\bowtie I $.
 Some references are
\cite{DF1, DF2, DFF1, DFF2, MY}. \\

 Let A and B be rings and let $\varphi: A\rightarrow B$ be a ring homomorphism making B an A-module. We say that A is a module retract of
 B if there exists a ring homomorphism $ \psi: B\rightarrow A$ such that $\psi o\varphi= id_{A}$. $\psi$
 is called  retraction of $\varphi$. See for instance \cite{G}. \\

 Our first main result in this paper is Theorem $2.1$ which gives
us a new characterization of $PF$-rings. Also, we provide a
necessary and sufficient condition for $R\bowtie I$ (resp., $R/I$
when $R$ is a Dedekind domain or $I$ is a primary ideal) to be
$PF$-ring. Our results generate new and original examples which
enrich the current literature with new families of $PF$-rings with zero-divisors.\\
 \end{section}

 \bigskip

 \begin{section}{Main Results}
\bigskip

 Recall that an R-module $M$ is called $P$-flat if, for any $(s,x) \in
R\times M$ such that $ sx=0$, then $ x\in$ $(0:s)M$. Now, we give
a new characterization for a class of $PF$-rings, which is the
first main result of this paper.

 \bigskip

  \begin{thm}
 Let $R$ be a commutative ring. Then the following conditions are equivalent:
\item [{\rm(1)}]Every ideal of $R$ is P-flat.
\item [{\rm(2)}]Every principal ideal of $R$ is P-flat.
\item [{\rm(3)}]$R$ is a $PF$-ring, that is every principal ideal of $R$ is flat.
\item [{\rm(4)}]For any elements $(s,x)\in R^{2}$ such that $sx=0$, there exists \\ $\alpha\in(0:s)$ such that  $x=\alpha
x$.
\end{thm}

\proof $(1)$ $\Longrightarrow$ (2) Clear.\\
$(2)$ $\Longrightarrow$ (3) Let $Ra$ be a principal ideal of $R$
generated by $a$. Our aim is to show that $Ra$ is flat. \\ Let $J$
be an ideal of $R$. We must show that
  $u:$ $Ra \otimes J \longrightarrow Ra\otimes R$,  where $u(a\otimes x)=ax$, is injective.
Let $a \in R$ and $x \in J$ such that $ax=0$. Hence, there exists
$\beta \in (0:x)$ and $ \lambda \in R$ such that $a=\beta\lambda
a$ (since $Ra$ is P-flat). Therefore,
$a\otimes x =\beta\lambda a\otimes x = \lambda a\otimes \beta x=0$, as desired.\\
$(3)$ $\Longrightarrow$ (4) Let $(s,x)$ be an element of $R^{2}$
such that $sx=0$. Our aim is to show that there exists $\beta\in
(0:s)$ such that $x=\beta x$.
 The principal ideal generated by $x$ is P-flat (since it is flat), so there exists $\alpha\in (0:s)$ and $r\in R$
 such that $x=\alpha rx =\beta x $ with $\beta =\alpha r \in (0:s)$.\\
$(4)$ $\Longrightarrow$ (1) Let I be an ideal of $R$. Let
$(s,x)\in R\times I$ such that $sx=0$. Hence, there exists
$\alpha\in(0:s)$ such that $ x=\alpha x$ and so $x\in (0:s)I$.
Therefore, I is P-flat, as desired. \qed

\bigskip

\begin{cor}
Let $R$ be a ring. The following conditions are equivalent:
\item [{\rm(1)}]  Every ideal of $R$ is $P$-flat.
\item [{\rm(2)}] Every ideal of $R_{Q}$ is $P$-flat for every prime ideal $Q$ of $R$.
\item [{\rm(3)}] Every ideal of $R_{m}$ is $P$-flat for every maximal ideal $m$ of $R$.
\item [{\rm(4)}] $R_{Q}$ is a domain for every prime ideal $Q$ of $R$.
\item [{\rm(5)}] $R_{m}$ is a domain for every maximal ideal $m$ of $R$.
\end{cor}

\proof By Theorem 2.1 and \cite[Theorem 4.2.2]{G}.\\
\qed

\bigskip

Recall that a ring $R$ is called an arithmetical ring if the
lattice formed by its ideals is distributive. If $wgldim(R) \leq
1$, then $R$ is an arithmetical ring. See for instance \cite{BKM, BG}. \\

Now, we add a condition with arithmetical in order to have
equivalence between arithmetical and $wgldim(R) \leq 1$.

\bigskip

\begin{prop} Let R be a ring. Then the following conditions are
equivalent:
\item [{\rm(1)}] $wgldim(R) \leq 1$.
\item [{\rm(2)}] $R$ is arithmetical and a $PF$-ring.
\item [{\rm(3)}] $R$ is arithmetical and every principal ideal of $R$ is flat.
\item [{\rm(4)}] $R$ is arithmetical and every principal ideal of $R$ is $P$-flat.
\item [{\rm(5)}] $R$ is arithmetical and every ideal of $R$ is $P$-flat.
\end{prop}

\proof $1) \Rightarrow 2) \Rightarrow 3) \Rightarrow 4)
\Rightarrow  5)$. By Theorem 2.1. \\
$5) \Rightarrow 1)$. Assume that the ring $R$ is arithmetical and
every ideal of $R$ is $P$-flat.  Our aim is to show that
$wgldim(R) \leq 1$. Let $I$ be a finitely generated ideal of $R$. Hence,
$I$ is $P$-flat and so $I$ is flat (since $R$ is arithmetical by
\cite[p. 236]{C}) and this completes the proof. \qed \\

\bigskip

 Now we show that the localization of a $PF$-ring is always a $PF$-ring.

\bigskip

\begin{prop}
 Let $R$ be a $PF$-ring and let S be a multiplicative subset of $R$. Then $S^{-1}(R)$ is a $PF$-ring.
\end{prop}

\proof  Assume that $R$ is a $PF$-ring and let $J$ be a principal
ideal of $S^{-1}(R)$. We claim that $J$ is flat. Indeed, since $J$
is a principal ideal of $S^{-1}R$, then there exists an element
$\dfrac{a}{b}$ of $J$ such that $J=S^{-1}(R)\dfrac{a}{b}$. Set $I
=Ra$. Hence, $I$ is flat since $R$ is a $PF$-ring and so $J
(=S^{-1}(I))$ is a flat ideal of $S^{-1}R$. It follows that
$S^{-1}(R)$ is a $PF$-ring. \qed
\\

\bigskip

Now, we study the transfer of $PF$-ring property to the direct
product.

\bigskip

\begin{prop}
Let $(R_{i})_{i\in I}$ be a family of commutative rings. Then
 $R=\prod_{i\in I}R_{i}$ is a $PF$-ring if and only if $R_{i}$ is
a $PF$-ring for all $i \in I$.
\end{prop}

\proof Assume that $R_i$ is a $PF$-ring for each $i \in I$ and set
$R =\prod_{i\in I}R_{i}$.  Let $x =(x_{i})_{i\in I}$ and $y
=(y_{i})_{i\in I}$ be two elements of $R$ such that $xy=0$. Then,
for every $i\in I,$ there exists $\alpha_{i}\in(0:x_{i})$ such
that $y_{i}=\alpha_{i}y_{i}$ (since $R_{i}$ is a $PF$-ring).
Hence, $(y_{i})_{i\in I}=(\alpha_{i})_{i\in I}(y_{i})_{i\in I}$
and $(\alpha_{i})_{i\in I}(x_{i})_{i\in I}=(\alpha_{i}x_{i})_{i\in
I}=0$. Therefore, $R$ is a $PF$-ring.\\

Conversely, assume that
$R =\prod_{i\in I}R_{i}$ is a $PF$-ring and we claim that $R_{i}$ is a $PF$-ring for every $i \in I$.\\
Indeed, let $i\in I $ and let $x_{i}$, $y_{i}$ be two elements of
$R_{i}$ such that $x_{i}y_{i}=0$. Consider $x=(a_{j})_{j\in I}$ ,
with $\left\{\begin{array}{c}
a_{i}=x_{i},\\
a_{j}=0  $ $ $ $ for $ $ j\neq i.\\
\end{array}\right. $
 and $y=(b_{j})_{j\in I}$ , with $\left\{\begin{array}{c}
b_{i}=y_{i}.\\
b_{j}=0 $ $ $ $ for $ $ i\neq j.\\\end{array}\right.$
 Since $R$ is a $PF$-ring, then there exists $\alpha\in (0:x)$ such that $y=\alpha y$ (that is, for all $j\in I$,  $b_{j}=\alpha_{j}b_{j}$ and
 $\alpha_{j}a_{j}=0$). Hence, $y_{i}=\alpha_{i}y_{i}$ with $\alpha_{i}\in (0:x_{i})$. Therefore, $R_{i}$ is a $PF$-ring for all $i\in I$
 and this completes the proof.\qed
\\

\bigskip

Next we study the transfer of $PF$-ring property to homomorphic
image. First, the following example shows that the homomorphic
image of a $PF$-ring is not always a $PF$-ring. \\

\bigskip

\begin{exam}
Let A be a domain and let $R=A[X]$. Then:
\item [{\rm(1)}] $R$ is a $PF$-ring since it is a domain.
\item [{\rm(2)}] $R/(X^{n})$ (for $n\geq2$) is not a $PF$-ring since
$\overline{X^{n}}=0 $ and $\overline{X}\neq0$.
\end{exam}

\bigskip

 Recall that if $R$ is a Dedekind domain and $I$ is a nonzero ideal of $R$, then $I=P_{1}^{\alpha_{1}}...P_{n}^{\alpha_{n}}$
 for some distinct prime ideals $P_{1},...,P_{n}$ uniquely determined by $I$ and
 some positive integers $\alpha_{1},...,\alpha_{n}$ uniquely determined
by $I$ (by \cite[Theorem 3.14]{G1}). \\

Now, when $R$ is a Dedekind domain or $I$ is a primary ideal, we
give a characterization of $R$ and $I$ such that $R/I$ is a $PF$-ring. \\

\bigskip

\begin{thm} Let $R$ be a ring and let $I$ be an ideal of $R$. Then:
 \item [{\rm(1)}] Assume that $R$ is a Dedekind domain and $I =P_{1}^{\alpha_{1}}...P_{n}^{\alpha_{n}}$ a nonzero ideal of $R$,
 where $P_{1},..,P_{n}$ are the prime ideals defined by $I$. Then $R/I$ is a $PF$-ring if and only if  $\alpha_{i}=1$ for all $i\in \{1,...,n\}$.
 \item [{\rm(2)}]  $I$ is a primary ideal of $R$ and $R/I$ is a $PF$-ring if and only if $I$ is a prime ideal of $R$.
  \end{thm}

\proof $1)$ Let $R$ be a Dedekind domain and let
$I=P_{1}^{\alpha_{1}}...P_{n}^{\alpha_{n}}$ for $P_{1},...,P_{n}$
be a nonzero prime ideals of $R$,
then  $R/I=\prod_{i=1}^n(R/P^{\alpha_{i}}_{i})$ .\\
 Assume that ${\alpha_{i}}=1$ for all $1\leq i\leq n$.
 Hence, $R/P_{i}$ is a $PF$-ring since $R/P_{i}$ is an integral domain, and so $ R/I=\prod_{i=1}^n(R/P_{i}^{\alpha_{i}})$ is a $PF$-ring by Proposition 2.5. \\

 Conversely, assume that $R/I=\prod_{i=1}^n(R/P^{\alpha_{i}}_{i})$ is a $PF$-ring. Let $i\in \{1,...,n\}$. Then
 $R/P_{i}^{\alpha_{i}}$  is a $PF$-ring by Proposition 2.5 . Hence, $R/P_{i}^{\alpha_{i}}$ is reduced
 and so the intersection of all prime
 ideals $Q$ of $R/P_{i}^{\alpha_{i}}$ is zero (i.e $\bigcap_{Q\in spect(R/P_{i}^{\alpha_{i}})}Q =\{0\}$) by \cite[Proposition 1.8]{A}.On the other hand for every prime ideals $Q$ of $R/P_{i}^{\alpha_{i}} $, there exist a prime ideal $Q'$ of $R$ such that $P_{i}^{\alpha_{i}}\subset Q'$ and $Q=Q'/P_{i}^{\alpha_{i}}$, then,
  $P_{i}/P_{i}^{\alpha_{i}}\subset Q$. It follows that
 $\{0\}=\bigcap_{Q\in spect(R/P_{i}^{\alpha_{i}})}Q $ =$P_{i}/P_{i}^{\alpha_{i}}$ and so
 $P_{i}=P_{i}^{\alpha_{i}}$, since $R$ is Dedekind domain then, $\alpha_{i}=1$. \\

$2)$ It's obvious that if $I$ is a prime ideal, then $R/I$ is
a $PF$-ring and $I$ is a primary ideal. \\
Conversely, assume that $I$ is a primary ideal and $R/I$ is a
$PF$-ring. Our aim is to show that $I$ is a prime ideal of $R$.
Let $x, y \in R$ such that $xy\in I$. We claim that $x\in I$ or
$y\in I$. Without loss of generality, we may assume that $x \notin
I$. Since $xy\in I$, then there exists an integer $n > 0$ such
that $y^{n}\in I$ (since $I$ is a primary ideal). Hence,
$\overline{y}^{n}=0$ and so $\overline{y}=0$ since $R/I$ is a
$PF$-ring; that is $y \in I$. Therefore, $x\in I$ or $y\in I$ and
so $I$ is a prime ideal of $R$, as desired.
 \qed
 \\

\bigskip

 Now, we are able to give examples of $PF$-rings and non-$PF$-rings. \\

\bigskip

 \begin{exam} (1) $\mathbb{Z}/4\mathbb{Z}$ is not a $PF$-ring by Theorem 2.8(1).\\
(2) $\mathbb{Z}/30\mathbb{Z}$ is a $PF$-ring by Theorem 2.8(1).
\end{exam}

\bigskip

 Now, we study the transfer of a $PF$-property to amalgamated duplication of a ring $R$ along an ideal $I$. \\

\bigskip

 \begin{thm}
 Let $R$ be a ring, and let $I$ be an ideal of $R$. Then,
   $R\bowtie I$ is a $PF$-ring if and only if $R$ is a $PF$ and $I$ is pure.
 \end{thm}

\bigskip

 We need the following lemma before proving this Theorem. \\

\bigskip

 \begin{lem} Let $R$ and $S$ be a rings and let $\varphi: R\rightarrow S$  be a ring homomorphism making $R$ a module retract of $S$.
 If $S$ is a $PF$-ring, then so is $R$.
 \end{lem}

 \proof  Let $\varphi: R\rightarrow S$ be a ring homomorphism and let $\psi: S\rightarrow R$ be a ring homomorphism such that
 $\psi o\varphi =id_{R}$. Let $(x,y)\in R^{2}$ such that $xy=0$. Then $\varphi (x) \varphi (y)= \varphi(xy)=0$ . Hence, there exists an element $\alpha \in S $ such that
 $\alpha\varphi(x)=0$ and $\varphi(y)=\alpha\varphi(y)$ (since $S$ is
 a $PF$-ring) and so $ y=\psi(\varphi(y))=\psi(\alpha\varphi(y))=\psi(\alpha)y$  and
 $\psi(\alpha)x=\psi(\alpha\varphi(x))=\psi(0)=0$, as desired. \qed\\

\bigskip

 \proof  of Theorem 2.9. \\
   Assume that $R\bowtie I$ is a $PF$-ring and we must to show that R is a $PF$-ring and $I$ is a pure ideal of $R$.
  We can easily show that $R$ is a module retract of $R\bowtie I$ where the retraction map $\varphi$
 is defined by $\varphi(r,r+i)=r$ and so $R$ is a $PF$-ring by Lemma 2.10. \\
 We claim that $I_{m}\in \{0,R_{m}\}$ for every maximal ideal $m$ of $R$. Let $m$ be an arbitrary maximal ideal of $R$, we have: $I\subseteq m$
 or $I\nsubseteq m$. If  $I\nsubseteq m$ then,  $I_{m}=R_{m}$. If $I\subseteq m$. Deny.  $I_{m} \notin
 \{0,R_{m}\}$ and so
 $(R\bowtie I)_{M}=R_{m}\bowtie I_{m}$, where $M$ a maximal ideal of $R\bowtie I$ such that $M\cap R=m$. Since $R_{m}$ is a domain,
 then $R_{m}\bowtie I_{m}$ is reduced and $O_{1} (=\{0\}\times I_{m})$ and $O_{2} (=I_{m}\times\{0\})$ are the only minimal prime ideals of
 $(R\bowtie I)_{M}$ by \cite[Proposition 2.1]{DF2}; hence it is not a $PF$-ring by \cite[Theorem 4.2.2]{G} (since $(R\bowtie I)_{M}$ is local), a
 desired contradiction. Therefore, $I_{m}\in \{0,R_{m}\}$ for every maximal ideal $m$ of $R$ . \\
 Conversely, assume that $R$ is a $PF$-ring and $I$ is a pure ideal of $R$, i.e. $I_{m}\in \{0,R_{m}\}$ for every maximal ideal $m$ of $R$ . Our aim is to prove
 that $R\bowtie I$ is a $PF$-ring. Using Corollary 2.2, we need to prove that $(R\bowtie I)_{M}$ is a $PF$-ring whenever $M$ is a maximal ideal
 of $R\bowtie I$. Let $M$ be an arbitrary maximal ideal of $R\bowtie I$ and set $m=M\cap R$. Then, necessarily $M\in\{M_{1},M_{2}\}$, where
 $M_{1}=\{(r,r+i)/r\in m,i\in I\}$ and $M_{2}=\{(r+i,r)/r\in m,r\in
 I\}$, by \cite[Theorem 3.5]{DF1}. On the other hand,
$I_{m}\in\{0,R_{m}\}$. Then,
 testing all cases of \cite[Proposition 7]{D}, we have two cases:\\
 (a) $(R\bowtie I)_{M}\cong R_{m}$  if $I_{m}=0$ or $I\nsubseteq m$.\\
 (b) $(R\bowtie I)_{M}\cong R_{m}\times R_{m}$  if $I_{m}=R_{m}$ and $I\subseteq m$.\\
 Since $R_{m}$ is a $PF$-ring (by Corollary 2.2), then so is $R_{m}\times R_{m}$ by Proposition 2.5 and hence $(R\bowtie I)_{M}$ is a $PF$-ring.\qed
\\

\bigskip

 \begin{cor}
 Let $R$ be a domain and let $I$ be a proper ideal of $R$. Then $R\bowtie I$ is never a $PF$-ring.
 \end{cor}

 \bigskip

 \begin{cor}
 Let $(R,m)$ be a local ring and let $I$ be a proper ideal of $R$. Then $R\bowtie I$ is never a $PF$-ring.
 \end{cor}

\bigskip

 Now we are able to construct a class of $PF$-rings. \\

\bigskip

 \begin{exam} Let $R$ be a $PF$-ring and let $I = Re$, where $e$ is an idempotent element of
 $R$. Then $R\bowtie I$ is a $PF$-ring by Theorem 2.9.
 \end{exam}

 \bigskip

 The following example shows that a subring of $PF$-ring is not always a $PF$-ring. For any ring
 $R$, we denote by $T(R)$ the total ring of quotients of $R$. \\

\bigskip

 \begin{exam} Let $R$ be an integral domain, $I$ a proper ideal of $R$  and let $S =R \bowtie
 I$. Then:
\item [{\rm(1)}] $S (=R\bowtie I)$ is not a $PF$-ring by Corollary
2.11.
 \item [{\rm(2)}]  $R\bowtie I \subseteq R \times R$ and $R \times
 R$ is a $PF$-ring by Proposition 2.5 (since $R$ is a $PF$-ring).
 \item [{\rm(3)}] $T(S) =T(R \times R) =K \times K$, where $K =T(R)$.
 \end{exam}

 \bigskip

 We end this paper by showing that the transfer of $PF$-ring property to Pullback is not always a
$PF$-ring.

\bigskip

 \begin{exam} Let $R$ be a domain and $I$ a proper ideal of $R$.
 Then:
 \item [{\rm(1)}]  The ring $R\bowtie I$ can be obtained as a pullback of $R$ and $R\times R$ over $
 R\times(R/I)$.
 \item [{\rm(2)}] The ring $R\bowtie I$ is not a $PF$-ring by Corollary 2.11.
 \item [{\rm(2)}] The rings $R$ and $R\times R$ are a $PF$-rings.
 \end{exam}

\end{section}
\bigskip
\bigskip




\bigskip\bigskip

\end{document}